\documentclass[twocolumn,aps,secnumarabic,nofootinbib,floatfix]{revtex4}
\usepackage{amssymb,amsmath2000,amsthm2000,times,mathptmx,pstricks,pst-node}
\newtheorem{theorem}{Theorem}[section]
\newtheorem{lemma}[theorem]{Lemma}

\newcommand{\C}{\mathbb{C}}
\renewcommand{\tensor}{\otimes}
\newcommand{\Inv}{\mathrm{Inv}}
\renewcommand{\sl}{\mathrm{sl}}
\newcommand{\so}{\mathrm{so}}
\renewcommand{\sp}{\mathrm{sp}}

\newcommand{\ie}{{\em i.e.}}

\newcommand{\evenodd}{\code{/fill {eofill} def /clip {eoclip} def}}

\let\blackandwhite\relax

\ifx\blackandwhite\undefined
  \psset{linecolor=blue}\newrgbcolor{darkred}{.5 0 0}
\else
  \psset{linecolor=black}\newgray{darkred}{0}
\fi

\psset{linewidth=.4pt,dash=3pt 3pt,doublesep=.05,dotsize=1pt 5}
\psset{hatchwidth=.3pt}
\SpecialCoor

\newcommand{\psgoesto}{\hspace{.5cm}\pspicture[.5](0,-.1)(1,.1)
\psline[linecolor=black]{->}(0,0)(1,0)
\endpspicture\hspace{.5cm}}
\newcommand{\rcrossing}{\pspicture[.4](-.6,-.5)(.6,.5)
\qline(.5;135)(.5;315)
\psline[border=.1](.5;45)(.5;225)
\endpspicture}
\newcommand{\lcrossing}{\pspicture[.4](-.6,-.5)(.6,.5)
\qline(.5;45)(.5;225)
\psline[border=.1](.5;135)(.5;315)
\endpspicture}

\newcommand{\vv}{\pspicture[.4](-.6,-.5)(.6,.5)
\psbezier(.5;45)(.25;45)(.25;315)(.5;315)
\psbezier(.5;225)(.25;225)(.25;135)(.5;135)
\endpspicture}
\newcommand{\hh}{\pspicture[.4](-.6,-.5)(.6,.5)
\psbezier(.5;45)(.25;45)(.25;135)(.5;135)
\psbezier(.5;225)(.25;225)(.25;315)(.5;315)
\endpspicture}
\newcommand{\rcurl}{\pspicture[.4](-.6,-.5)(.6,.5)
\psbezier(.5;120)(-.25,0)(.25,-.5)(.25,0)
\psbezier[border=.1](.5;240)(-.25,0)(.25,.5)(.25,0)
\endpspicture}
\newcommand{\leftv}{\pspicture[.4](-.6,-.5)(.6,.5)
\psbezier(.5;240)(.25;240)(.25;120)(.5;120)
\endpspicture}
\newcommand{\singleloop}{\pspicture[.4](-.6,-.5)(.6,.5)
\pscircle(0,0){.4}
\endpspicture}
\newcommand{\doubleloop}{\pspicture[.4](-.6,-.5)(.6,.5)
\pscircle[doubleline=true](0,0){.4}
\endpspicture}

\newcommand{\doublesquare}{\pspicture[.4](-.9,-.9)(.9,.9)
\psline[doubleline=true](.4;45)(.8;45)
\psline[doubleline=true](.4;135)(.8;135)
\psline[doubleline=true](.4;225)(.8;225)
\psline[doubleline=true](.4;315)(.8;315)
\psline(.4;45)(.4;135)\psline(.4;135)(.4;225)
\psline(.4;225)(.4;315)\psline(.4;315)(.4;45)
\endpspicture}
\newcommand{\hvert}{\pspicture[.4](-.6,-.7)(.6,.7)
\qline(0, .25)(.35, .6)\qline(0, .25)(-.35, .6)
\qline(0,-.25)(.35,-.6)\qline(0,-.25)(-.35,-.6)
\psline[doubleline=true](0,-.25)(0,.25)
\endpspicture}
\newcommand{\hhoriz}{\pspicture[.4](-.8,-.5)(.8,.5)
\qline( .25,0)( .6,.35)\qline( .25,0)( .6,-.35)
\qline(-.25,0)(-.6,.35)\qline(-.25,0)(-.6,-.35)
\psline[doubleline=true](-.25,0)(.25,0)
\endpspicture}
\newcommand{\bwhvert}{\pspicture[.4](-.6,-.8)(.6,.8)
\pscustom[fillstyle=crosshatch,linestyle=none]{
\psarc(0,-.25){.5}{225}{315}
\psline(.35,-.6)(0,-.25)(-.35,-.6)}
\psline(.35,-.6)(0,-.25)(-.35,-.6)
\pscustom[fillstyle=crosshatch,linestyle=none]{
\psarc(0,.25){.5}{45}{135}
\psline(-.35,.6)(0,.25)(.35,.6)}
\psline(-.35,.6)(0,.25)(.35,.6)
\psline[doubleline=true](0,-.25)(0,.25)
\endpspicture}
\newcommand{\bwhhoriz}{\pspicture[.4](-.9,-.7)(.9,.7)
\pscustom[fillstyle=crosshatch,linestyle=none]{
\psline(.6,-.35)(.25,0)(.6,.35)
\psarc(0,0){.695}{30.256}{149.744}
\psline(-.6,.35)(-.25,0)(-.6,-.35)
\psarc(0,0){.695}{-149.744}{-30.256}}
\psline(-.6,-.35)(-.25,0)(-.6,.35)
\psline(.6,.35)(.25,0)(.6,-.35)
\psline[doubleline=true](-.25,0)(.25,0)
\endpspicture}
\newcommand{\bwhh}{\pspicture[.4](-.6,-.5)(.6,.5)
\pscustom[fillstyle=crosshatch,linestyle=none]{
\psbezier(.5;135)(.25;135)(.25;45)(.5;45)
\psarc(0,0){.5}{45}{135}
\psbezier[liftpen=2](.5;315)(.25;315)(.25;225)(.5;225)
\psarc(0,0){.5}{225}{315}}
\psbezier(.5;45)(.25;45)(.25;135)(.5;135)
\psbezier(.5;225)(.25;225)(.25;315)(.5;315)
\endpspicture}
\newcommand{\bwvv}{\pspicture[.4](-.6,-.5)(.6,.5)
\pscustom[fillstyle=crosshatch,linestyle=none]{
\psbezier(.5;315)(.25;315)(.25;45)(.5;45)
\psarc(0,0){.5}{45}{135}
\psbezier(.5;135)(.25;135)(.25;225)(.5;225)
\psarc(0,0){.5}{225}{315}}
\psbezier(.5;45)(.25;45)(.25;315)(.5;315)
\psbezier(.5;225)(.25;225)(.25;135)(.5;135)
\endpspicture}

\begin{document}

\title{Jaeger's Higman-Sims state model and the $B_2$ spider}
\author{Greg Kuperberg}
\thanks{The author was supported by NSF grant \#DMS-9423300.}
\affiliation{Department of Mathematics, University of California,
     Davis, CA 95616}
\email[Email:]{greg@math.ucdavis.edu}

\begin{abstract}
Jaeger \cite{Jaeger:spin} discovered a remarkable checkerboard state
model based on the Higman-Sims graph that yields a value of the
Kauffman polynomial, which is a quantum invariant of links.  We
present a simple argument that the state model has the desired
properties using the combinatorial $B_2$ spider \cite{Kuperberg:spiders}.
\end{abstract}
\maketitle

\section{Introduction}

Two related approaches to defining quantum topological invariants
are skein relations and state models.  One important example of the
former is the Kauffman polynomial, while an important class of 
the latter is the class of checkerboard state models.

Given two numbers or indeterminates $Q$ and $d$, the Kauffman polynomial
is a function on link projections on the 2-sphere defined axiomatically
by the rules:
\begin{align}
\rcrossing - \lcrossing &= (Q-Q^{-1})\left(\vv - \hh\right) \label{eskein} \\
\singleloop &= \frac{Q^{d-1} + Q - Q^{-1} - Q^{1-d}}{Q - Q^{-1}} \nonumber\\
\rcurl &= Q^{d-1}\leftv \nonumber
\end{align}
by the rule that it is invariant under the second and third Reidemeister
moves:
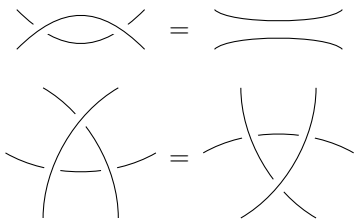
\begin{align*}
\pspicture[.4](-1.1,-.5)(1.1,.5)
\pccurve[angleA=-45,angleB=-135,ncurv=1](-.85,.264)(.85,.264)
\pccurve[border=.1,angleA=45,angleB=135,ncurv=1](-.85,-.264)(.85,-.264)
\endpspicture &=
\pspicture[.4](-1.1,-.5)(1.1,.5)
\pccurve[angleA=-45,angleB=-135,ncurv=.33](-.85,.264)(.85,.264)
\pccurve[border=.1,angleA=45,angleB=135,ncurv=.33](-.85,-.264)(.85,-.264)
\endpspicture \\
\pspicture[.42](-1.1,-1)(1.1,1)
\pcarc[border=.1,arcangle=30](1;0)(1;180)
\pcarc[border=.1,arcangle=30](1;120)(1;300)
\pcarc[border=.1,arcangle=30](1;240)(1;60)
\endpspicture &=
\pspicture[.42](-1.1,-1)(1.1,1)
\pcarc[border=.1,arcangle=30](1;180)(1;0)
\pcarc[border=.1,arcangle=30](1;300)(1;120)
\pcarc[border=.1,arcangle=30](1;60)(1;240)
\endpspicture
\end{align*}
and by the rule that its value at the empty link is 1.   These rules are an
example of a skein theory, a concept which can be understood with some
elementary background.

A (tame) knot or link is represented by a knot projection (a tetravalent graph
embedded in the 2-sphere with vertices decorated to distinguish over-crossings
from under-crossings); it is known that a function on links which is invariant
under the three Reidemeister moves (indicated above) is a function on  links. 
It is further understood that invariance under the second and third
Reidemeister moves (regular isotopy invariance) is almost as strong as
invariance under all three, in a manner analogous to the difference between
linear and projective representations of a group.  In this paper, we will
loosely call a regular isotopy invariant a link invariant.  A skein theory
describes a function on knot projections by axioms relating projections that
differ only in a small region.  Typically one implicitly considers an invariant
$I$ and one writes graphical equations where a knot projection $P$ denotes
$I(P)$.  Furthermore, if an equation involves link fragments (tangles), the
fragments should all have the same boundary, and the equation is  read as
relating link projections that differ only in the indicated fragments.  For
example, equation~\ref{eskein} is really infinitely many equations relating any
quadruple of link projections that differ in only one crossing.

Kauffman \cite{Kauffman:regular} proved that the Kauffman polynomial exists
uniquely; \ie, that the skein relations are consistent and complete. It is
therefore a Laurent polynomial in $Q$ and $Q^d$. (Our parameterization of the
Kauffman polynomial is slightly different from Kauffman's.)  The specialization
$d = -2$ is called the Kauffman bracket, another polynomial which is up to
normalization the same as the Jones polynomial.  (Both the Kauffman polynomial
and the Kauffman bracket are clearly not invariant under the first Reidemeister
move, but a simple extra normalization factor achieves this invariance as well
and produces a function fully invariant under link isotopy.  This is the main
difference between the Kauffman bracket and the Jones polynomial.)

An invariant with a skein theory often has an alternative definition using a
state model; we will consider a particular type of state model for link
projections called a checkerboard model\cite{Jones:pacific}.  A checkerboard
model is given by a state set $S$, a number $x$, and two symmetric functions
$W_+$ and $W_-$ from $S \times S$ to a commutative ring $R$.  Given a link
projection on the 2-sphere, a checkerboard coloring is one of the two
alternating black-white colorings of the complementary regions:
$$
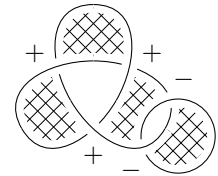
\pspicture(-1.2,-1.4)(1.2,1.6)
\pscustom[fillstyle=crosshatch,linestyle=none]{\evenodd
\psbezier(1.1;90)(1.27;60)(.7;30)(.35;330)
\psbezier(.35;330)(.7;270)(1.27;240)(1.1;210)
\psbezier(1.1;210)(1.27;180)(.7;150)(.35;90)
\psbezier(.35;90)(.7;30)(1.27;0)(1.1;330)
\psbezier(1.1;330)(1.27;300)(.7;270)(.35;210)
\psbezier(.35;210)(.7;150)(1.27;120)(1.1;90)
\psarc[liftpen=2](1.3;330){.5}{0}{360}}
\psbezier[border=.1](.35;330)(.7;270)(1.27;240)(1.1;210)
\psbezier[border=.1](.35;90)(.7;30)(1.27;0)(1.1;330)
\psbezier[border=.1](.35;210)(.7;150)(1.27;120)(1.1;90)
\pscircle[border=.1](1.3;330){.5}
\psbezier[border=.1](1.1;90)(1.27;60)(.7;30)(.35;330)
\psbezier[border=.1](1.1;210)(1.27;180)(.7;150)(.35;90)
\psbezier[border=.1](1.1;330)(1.27;300)(.7;270)(.35;210)
\rput(.5,-1.1){$-$}\rput(1.2,.1){$-$}
\rput(.9;270){$+$}\rput(.9;150){$+$}\rput(.9;30){$+$}
\endpspicture
$$
Given a checkerboard coloring, the crossings can be labelled as positive or
negative, as indicated, according to their sense relative to the neighboring
black regions.  (The labelling depends on an orientation of the 2-sphere.) The
black regions are called the atoms of the state model.  A state is a function
from the atoms to the state set.  Given a state, the weight of a positive
(respectively negative) crossing, which takes values in $\C$ or some other
field, is given by some function $W_+(a,b)$ (resp. $W_-(a,b)$) if the two
incident atoms are assigned states $a$ and $b$; these functions are called
interactions. The weight of a state is then the product of the weights of the
atoms, and the state sum $Z$ is the total weight of all states.  If $\chi$ is
the total Euler characteristic of all black regions, it may happen that the
normalized state sum $Z' = x^{-\chi} Z$ is a regular isotopy invariant and in
particular one that satisfies skein relations.

For example, the Potts model is a checkerboard model whose normalized
state sum is a value of the Kauffman bracket, and in particular is a regular
isotopy invariant.  Choose a real or complex $q$ such that $n =
q+2+q^{-1}$ is a positive integer and choose a root $q^{1/4}$. Then the Potts
model of order $n$ has a state set with $n$ elements and the following
weights:
\begin{align*}
W_+(a,a) &= q^{3/4} \\
W_-(a,a) &= q^{-3/4} \\
W_+(a,b) &= -q^{-1/4} \qquad \mbox{when $a \ne b$} \\
W_-(a,b) &= -q^{1/4} \qquad \mbox{when $a \ne b$}
\end{align*}
The Euler normalization $x = -(q^{1/2} + q^{-1/2}) = \pm \sqrt{n}$.
One can check that the normalized state sum is then
the Kauffman bracket with $Q = q^{1/4}$.

There are only a few known non-trivial checkerboard models that produce
link invariants \cite{dHJ:graph}. By far the most interesting of these
is the Higman-Sims state model discovered by Jaeger \cite{Jaeger:spin}.  This
model produces the value of the Kauffman polynomial at $Q = \tau$, the golden
ratio, and $d = -4$.  In quantum groups terms, this corresponds to the
$q = \tau^2$ point of the quantum group $U_q(\sp(4))$.

In this paper, we present an alternative argument that the Higman-Sims state
model produces a topological invariant and a value of the Kauffman polynomial. 
The idea is to use the combinatorial $B_2$ spider \cite{Kuperberg:spiders}, a
skein theory of graphs that produces the Kauffman polynomial for arbitrary $Q$
with $d = -4$.  The author hopes that this argument will help draw attention to
Jaeger's remarkable state model and help elucidate a mysterious relationship
between the Higman-Sims sporadic simple group and the quantum group
$U_q(\sp(4))$.  We also refer the reader to another review of Jaeger's model by
de la Harpe \cite{delaHarpe:jaeger}.

\section{Spiders, skein modules, and invariants}

The combinatorial $B_2$ spider is essentially a skein theory of trivalent,
planar graphs (called webs) with unoriented edges of two types, which are
called type 1 and type 2 strands and are denoted by single and double edges. 
Such  a skein theory uses the same concepts and allows the same notation as a
skein theory for link projections, except that it describes a function on a
class of planar graphs and we require no particular topological invariance a
priori (other than invariance under isotopy in the 2-sphere). The only allowed
vertices in $B_2$ webs are those with two single edges and one double edge:
$$
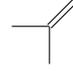
\pspicture(-.7,-.7)(.7,.7)
\qline(-.5,0)(0,0)\qline(0,-.5)(0,0)\psline[doubleline=true](0,0)(.35,.35)
\endpspicture
$$
The skein relations are:
\begin{align}
\singleloop &= -(q^2+q+q^{-1}+q^{-2}) \label{eloop} \\
\doubleloop &= q^3+q+1+q^{-1}+q^{-3} \nonumber \\
\pspicture[.4](-.6,-.5)(.6,.5)
\psline[doubleline=true](-.5,0)(0,0)
\psbezier(0,0)(.7,.7)(.7,-.7)(0,0)
\endpspicture
&= 0 \nonumber \\
\pspicture[.4](-.8,-.5)(.8,.5)
\psline[doubleline=true](-.7,0)(-.3,0)
\pcarc[arcangle=45](-.3,0)(.3,0)
\pcarc[arcangle=-45](-.3,0)(.3,0)
\psline[doubleline=true](.3,0)(.7,0)
\endpspicture
&= -(q+2+q^{-1})\pspicture[.4](-.6,-.5)(.6,.5)
\psline[doubleline=true](-.5,0)(.5,0)
\endpspicture \nonumber \\
\pspicture[.4](-.9,-.9)(.9,.9)
\psline[doubleline=true](.4;90)(.8;90)
\psline[doubleline=true](.4;210)(.8;210)
\psline[doubleline=true](.4;330)(.8;330)
\pcarc[arcangle=-15](.4;90)(.4;210)
\pcarc[arcangle=-15](.4;210)(.4;330)
\pcarc[arcangle=-15](.4;330)(.4;90)
\endpspicture
&= 0 \nonumber \\
\hvert &- \hhoriz = \hh - \vv \label{eexchange}
\end{align}
Here $q$ is a complex number or an indeterminate with a preferred square root
$q^{1/2}$.  Together with the stipulation that the empty web has value 1, these
skein relations again define a unique function on $B_2$ webs
\cite{Kuperberg:g2}.

To understand something of the relation between the $B_2$ spider and the Lie
algebra $B_2$, we can consider skein modules, which are another general concept
in skein theory.  Instead of considering functions on link projections or $B_2$
webs, we consider formal linear combinations of such objects over a field such
as $\C(q^{1/2})$, or sometimes over a ring.  We quotient the space of all
formal linear combinations by the skein relations, more properly interpreted as
relators.  Discarding stipulations about empty diagrams, we can say that the
skein module for the Kauffman polynomial and the skein module for the $B_2$
spider are both 1-dimensional.  More generally, if we fix a boundary of a
tangle or a $B_2$ web (meaning that the tangle or web is embedded in a disk and
has prespecified univalent endpoints on the boundary of the disk), we can
consider the skein module of tangles or webs with this boundary.

In the $B_2$ spider, the skein module is called a web space and all of its
elements are called webs.  If a web space has $n$ endpoints of type 1 and $k$
of type 2, then it is isomorphic to the invariant space
$\Inv(V(\lambda_1)^{\tensor k} \tensor V(\lambda_2)^{\tensor n})$, where
$V(\lambda_1)$ and $V(\lambda_2)$ are, respectively, the 5-dimensional and
4-dimensional irreducible representations of the quantum group $U_q(\sp(4))$
\cite{Kuperberg:spiders}.  (Such a quantum group is an algebra which
specializes to the usual universal enveloping algebra when $q=1$.)  Using
$U_q(\sp(4))$, one can define an algebraic $B_2$ spider, and one can say that
the algebraic and combinatorial $B_2$ spiders are isomorphic when $q$ is a
transcendental element in a field or is not a root of unity.

Working in the $B_2$ spider, we may define particular webs called crossings by
the equations
\begin{align}
\rcrossing =&\ -q^{1/2}\vv - \frac{q^{-1}}{q^{1/2} + q^{-1/2}}\hh \nonumber \\
&\ - \frac{1}{q^{1/2} + q^{-1/2}}\hvert \label{ecrossing} \\
\pspicture[.4](-.6,-.5)(.6,.5)
\psline[doubleline=true](.5;135)(.5;315)
\psline[doubleline=true,border=.1](.5;45)(.5;225)
\endpspicture
=&\ q
\pspicture[.4](-.6,-.5)(.6,.5)
\psbezier[doubleline=true](.5;45)(.25;45)(.25;315)(.5;315)
\psbezier[doubleline=true](.5;225)(.25;225)(.25;135)(.5;135)
\endpspicture
+ q^{-1}
\pspicture[.4](-.6,-.5)(.6,.5)
\psbezier[doubleline=true](.5;45)(.25;45)(.25;135)(.5;135)
\psbezier[doubleline=true](.5;225)(.25;225)(.25;315)(.5;315)
\endpspicture \nonumber \\
&\ + \frac{1}{q + 2 +q^{-1}} \doublesquare \nonumber
\end{align}
These webs then satisfy the regular isotopy equations and yield invariants of
framed graphs in $S^3$ or in a ball which take values in the appropriate web
spaces. Moreover, a crossing of two type 1 strands satisfies Kauffman's
relations with  $Q = q^{1/2}$ and $d = -4$, while a crossing of two type 2
strands satisfies Kauffman's relation with $Q = q$ and $d = 5$.  So one can say
that there is a homomorphism (in fact an epimorphism which is far from
injective) from a one-variable slice of the Kauffman spider to the $B_2$
spider.

\section{The model}

Given a $B_2$ web on the sphere or in a disk, a checkerboard coloring is a
coloring of its faces such that two regions that meet at a type 1 strand have
opposite colors, while two regions that meet at a type 2 strand have the same
color.  For technical reasons we only consider webs with no type 2 strands
at the boundary and with no closed loops of type 2.  We can then write the
same $B_2$ skein relations for this class of colored webs, for example:
\begin{align*}
\pspicture[.4](-1.3,-.5)(1.3,.5)
\pscustom[fillstyle=crosshatch,linestyle=none]{
\psarc(-.7,0){.5}{135}{225}
\psline(-1.05,-.35)(-.7,0)(-1.05,.35)}
\psline(-1.05,-.35)(-.7,0)(-1.05,.35)
\pscustom[fillstyle=crosshatch,linestyle=none]{
\psarc(.7,0){.5}{315}{45}
\psline(1.05,.35)(.7,0)(1.05,-.35)}
\psline(1.05,.35)(.7,0)(1.05,-.35)
\psline[doubleline=true](-.7,0)(-.3,0)
\pscustom[fillstyle=crosshatch]{
\psbezier(-.3,0)(-.15,.26)(.15,.26)(.3,0)
\psbezier(.3,0)(.15,-.26)(-.15,-.26)(-.3,0)}
\psline[doubleline=true](.3,0)(.7,0)
\psline(1.05,-.35)(.7,0)(1.05,.35)
\endpspicture
&= -(q+2+q^{-1})
\pspicture[.4](-.9,-.5)(.9,.5)
\pscustom[fillstyle=crosshatch,linestyle=none]{
\psarc(-.25,0){.5}{135}{225}
\psline(-.6,-.35)(-.25,0)(-.6,.35)}
\psline(-.6,-.35)(-.25,0)(-.6,.35)
\pscustom[fillstyle=crosshatch,linestyle=none]{
\psarc[linewidth=0](.25,0){.5}{315}{45}
\psline(.6,.35)(.25,0)(.6,-.35)}
\psline(.6,.35)(.25,0)(.6,-.35)
\psline[doubleline=true](-.25,0)(.25,0)
\endpspicture \\
\pspicture[.4](-1.2,-1.2)(1.2,1.2)
\pscustom[fillstyle=crosshatch,linestyle=none]{\evenodd
\psarc(0,0){1.15}{345}{75}
\psline(1.15;75)(.8;90)(1.15;105)
\psarc(0,0){1.15}{105}{195}
\psline(1.15;195)(.8;210)(1.15;225)
\psarc(0,0){1.15}{225}{315}
\psline(1.15;315)(.8;330)(1.15;345)
\psbezier[liftpen=2](.4;90)(.3;130)(.3;170)(.4;210)
\psbezier(.4;210)(.3;250)(.3;290)(.4;330)
\psbezier(.4;330)(.3;10)(.3;50)(.4;90)}
\psline(1.15;75)(.8;90)(1.15;105)
\psline(1.15;195)(.8;210)(1.15;225)
\psline(1.15;315)(.8;330)(1.15;345)
\psline[doubleline=true](.4;90)(.8;90)
\psline[doubleline=true](.4;210)(.8;210)
\psline[doubleline=true](.4;330)(.8;330)
\psbezier(.4;90)(.3;130)(.3;170)(.4;210)
\psbezier(.4;210)(.3;250)(.3;290)(.4;330)
\psbezier(.4;330)(.3;10)(.3;50)(.4;90)
\endpspicture
&= 0
\end{align*}
These relations constitute a perfectly valid skein theory even though
the colorings of the faces are somewhat redundant.

We consider state models on checkerboard colorings of $B_2$ webs. Let $S$ be
a state set; as before, a state is a function from the black regions (atoms)
to the state set. Assuming that there are no crossings, we consider two
interactions (functions) $W_H$ and $W_I$ from $S \times S$ to $\C$.  The weight of a
state has a factor of $W_I(a,b)$ for every type 2 edge that bridges an atom
in state $a$ with an atom in state $b$ and a factor of $W_H(a,b)$ for every
type 2 edge that lies on the border between an atom in state $a$ from an atom
in state $b$:
$$
\pspicture(-.6,-1.3)(.6,.8)
\pscustom[fillstyle=crosshatch,linestyle=none]{
\psarc(0,-.25){.5}{225}{315}
\psline(.35,-.6)(0,-.25)(-.35,-.6)}
\psline(.35,-.6)(0,-.25)(-.35,-.6)
\pscustom[fillstyle=crosshatch,linestyle=none]{
\psarc(0,.25){.5}{45}{135}
\psline(-.35,.6)(0,.25)(.35,.6)}
\psline(-.35,.6)(0,.25)(.35,.6)
\psline[doubleline=true](0,-.25)(0,.25)
\cput*(0,.7){$a$}
\cput*(0,-.7){$b$}
\rput(0,-1.3){$W_I(a,b)$}
\endpspicture
\qquad
\pspicture(-.9,-1.3)(.9,.8)
\pscustom[fillstyle=crosshatch,linestyle=none]{
\psline(.6,-.35)(.25,0)(.6,.35)
\psarc(0,0){.695}{30.256}{149.744}
\psline(-.6,.35)(-.25,0)(-.6,-.35)
\psarc(0,0){.695}{-149.744}{-30.256}}
\psline(-.6,-.35)(-.25,0)(-.6,.35)
\psline(.6,.35)(.25,0)(.6,-.35)
\psline[doubleline=true](-.25,0)(.25,0)
\cput*(0,.4){$a$}
\cput*(0,-.4){$b$}
\rput(0,-1.3){$W_H(a,b)$}
\endpspicture
$$
As before, the state sum $Z$ is the total weight of all states and $Z' =
x^{-\chi} Z$ is the normalized state sum for some constant $x$.  Note that for
a web in a disk, for every fixed state of all atoms at the boundary, there is a
state sum over all states in the interior; we do not sum over the colorings of
the boundary regions.  If we establish that the state sum $Z'$ satisfies the
above skein relations (equations~(\ref{eloop}) to (\ref{eexchange})), then we
can say that it is a linear functional on web spaces. We can further define
interactions $W_-$ and $W_+$ for crossings by using equation~(\ref{ecrossing})
(see also the relation between $W_H$ and $W_I$ below):
$$
W_{\pm} = -q^{\pm 1/2} x I - \frac{q^{\mp 1}}{q^{1/2}+q^{-1/2}}N -
\frac{W_I}{q^{1/2}+q^{-1/2}}
$$
Here and below the interactions are interpreted as $n \times n$ matrices, $I$
is the identity matrix, and $N$ is the matrix of all 1's. The interactions
$W_\pm$ then constitute a checkerboard model for link projections, one whose
normalized state sum is automatically a value of the Kauffman polynomial at
$d=-4$.

Consider the restrictions on $W_H$, $W_I$, and $x$ given by the skein
relations. Assume that the state set $S$ has $n$ elements.  Let
\begin{align*}
h &= -(q+2+q^{-1}) \\
\ell &= -(q^2+q+q^{-1}+q^{-2})
\end{align*}
The relations
$$
\pspicture[.4](-.6,-.5)(.6,.5)
\pscircle[fillstyle=crosshatch](0,0){.4}
\endpspicture
= \ell
\hspace{2cm}
\pspicture[.4](-.8,-.7)(.8,.7)
\pscustom[fillstyle=crosshatch,linestyle=none]{\evenodd
\psframe(-.7,-.7)(.7,.7)
\psarc[liftpen=2](0,0){.4}{0}{360}}
\pscircle(0,0){.4}
\endpspicture
= \ell
$$
say that $\ell = \frac{n}{x}$ and $\ell = x$, which implies that $n = \ell^2$.
The relation
$$
\pspicture[.4](-.6,-.5)(.8,.5)
\pscustom[fillstyle=crosshatch,linestyle=none]{\evenodd
\psframe(-.5,-.5)(.7,.5)
\psbezier(0,0)(.7,.7)(.7,-.7)(0,0)}
\psline[doubleline=true](-.5,0)(0,0)
\psbezier(0,0)(.7,.7)(.7,-.7)(0,0)
\endpspicture
= 0
$$
says that $W_H(a,a) = 0$.  The relation 
$$
\pspicture[.4](-1.3,-1)(1.3,1)
\pscustom[fillstyle=crosshatch,linestyle=none]{\evenodd
\psline(-1.05,-.35)(-.7,0)(-1.05,.35)
\psbezier(-1.05,.35)(-.35,1.05)(.35,1.05)(1.05,.35)
\psline(1.05,.35)(.7,0)(1.05,-.35)
\psbezier(1.05,-.35)(.35,-1.05)(-.35,-1.05)(-1.05,-.35)
\psbezier[liftpen=2](-.3,0)(-.15,.26)(.15,.26)(.3,0)
\psbezier(.3,0)(.15,-.26)(-.15,-.26)(-.3,0)}
\psline(-1.05,.35)(-.7,0)(-1.05,-.35)
\psline[doubleline=true](-.7,0)(-.3,0)
\psbezier(-.3,0)(-.15,.26)(.15,.26)(.3,0)
\psbezier(.3,0)(.15,-.26)(-.15,-.26)(-.3,0)
\psline[doubleline=true](.3,0)(.7,0)
\psline(1.05,-.35)(.7,0)(1.05,.35)
\endpspicture
 = h \bwhhoriz
$$
says that $W_H(a,b)$ is either 0 or $h$ for all $a$ and $b$.  Since $W_H(a,b)$
is symmetric, it is proportional to the adjacency matrix $A$ of some graph
$J$ with vertex set $S$.  (The graph $J$ need not be planar or otherwise
resemble a web.)  The relation
$$
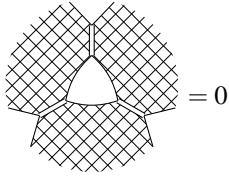
\pspicture[.4](-1.2,-1.2)(1.2,1.2)
\pscustom[fillstyle=crosshatch,linestyle=none]{\evenodd
\psarc(0,0){1.15}{345}{75}
\psline(1.15;75)(.8;90)(1.15;105)
\psarc(0,0){1.15}{105}{195}
\psline(1.15;195)(.8;210)(1.15;225)
\psarc(0,0){1.15}{225}{315}
\psline(1.15;315)(.8;330)(1.15;345)
\psbezier[liftpen=2](.4;90)(.3;130)(.3;170)(.4;210)
\psbezier(.4;210)(.3;250)(.3;290)(.4;330)
\psbezier(.4;330)(.3;10)(.3;50)(.4;90)}
\psline(1.15;75)(.8;90)(1.15;105)
\psline(1.15;195)(.8;210)(1.15;225)
\psline(1.15;315)(.8;330)(1.15;345)
\psline[doubleline=true](.4;90)(.8;90)
\psline[doubleline=true](.4;210)(.8;210)
\psline[doubleline=true](.4;330)(.8;330)
\psbezier(.4;90)(.3;130)(.3;170)(.4;210)
\psbezier(.4;210)(.3;250)(.3;290)(.4;330)
\psbezier(.4;330)(.3;10)(.3;50)(.4;90)
\endpspicture = 0
$$
says that $J$ is triangle-free.  The relation
$$ \bwhvert - \bwhhoriz = \bwhh - \bwvv $$
reads algebraically as
$$ W_I - W_H = N - xI.$$
where $W_I$ and $W_H$ can be interpreted as matrices, $N$ is the matrix whose
entries are all 1, and $I$ is the identity matrix.  This equation can be
taken as a definition of $W_I$ in terms of $W_H = hA$:
\begin{equation}
W_I = hA + N - xI \label{ewi}
\end{equation}
The relation
$$
\pspicture[.4](-.6,-.5)(.6,.5)
\psline[doubleline=true](-.5,0)(0,0)
\psbezier[fillstyle=crosshatch](0,0)(.7,.7)(.7,-.7)(0,0)
\endpspicture = 0
$$
is equivalent to the equation $W_IN = 0$.  Using equation~(\ref{ewi})
and the identity $N^2 = nN$, we obtain
$$AN = \frac{n-x}{h}N.$$
This equation says that $J$ is regular (1-point regular as defined
below) and the degree of a vertex is $v = \frac{n-x}{h}$.  The relation
$$
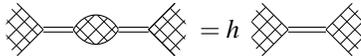
\pspicture[.4](-1.3,-.5)(1.3,.5)
\pscustom[fillstyle=crosshatch,linestyle=none]{
\psarc(-.7,0){.5}{135}{225}
\psline(-1.05,-.35)(-.7,0)(-1.05,.35)}
\psline(-1.05,-.35)(-.7,0)(-1.05,.35)
\pscustom[fillstyle=crosshatch,linestyle=none]{
\psarc(.7,0){.5}{315}{45}
\psline(1.05,.35)(.7,0)(1.05,-.35)}
\psline(1.05,.35)(.7,0)(1.05,-.35)
\psline[doubleline=true](-.7,0)(-.3,0)
\pscustom[fillstyle=crosshatch]{
\psbezier(-.3,0)(-.15,.26)(.15,.26)(.3,0)
\psbezier(.3,0)(.15,-.26)(-.15,-.26)(-.3,0)}
\psline[doubleline=true](.3,0)(.7,0)
\psline(1.05,-.35)(.7,0)(1.05,.35)
\endpspicture
 = h
\pspicture[.4](-.9,-.5)(.9,.5)
\pscustom[fillstyle=crosshatch,linestyle=none]{
\psarc(-.25,0){.5}{135}{225}
\psline(-.6,-.35)(-.25,0)(-.6,.35)}
\psline(-.6,-.35)(-.25,0)(-.6,.35)
\pscustom[fillstyle=crosshatch,linestyle=none]{
\psarc[linewidth=0](.25,0){.5}{315}{45}
\psline(.6,.35)(.25,0)(.6,-.35)}
\psline(.6,.35)(.25,0)(.6,-.35)
\psline[doubleline=true](-.25,0)(.25,0)
\endpspicture
$$
reads algebraically as the matrix equation
$$W_I^2 = xh W_I,$$
which implies that $W_I$ as a linear operator has only two eigenvalues, namely
0 and $xh$. Moreover, given that $J$ is regular, the property that some
linear combination of $A$, $N$, and $I$ has a quadratic minimal polynomial is
equivalent to the property that $J$ is strongly regular \cite{Jaeger:spin}, or
2-point regular as defined below.  Thus, we see that the $B_2$ skein relations
determine the parameters $x$ and $n$ of a checkerboard state model in terms
of $q$, that they imply that the model is essentially determined by a certain
graph $J$, and that they place strong restrictions on $J$.

Finally, there is the relation
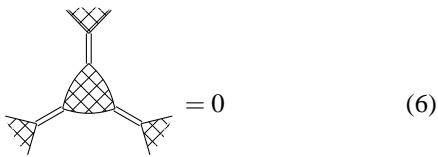
\begin{equation}
\pspicture[.4](-1.2,-1.2)(1.2,1.2)
\pscustom[fillstyle=crosshatch,linestyle=none]{
\psarc(0,0){1.15}{75}{105}
\psline(1.15;75)(.8;90)(1.15;105)
\psarc[liftpen=2](0,0){1.15}{195}{225}
\psline(1.15;195)(.8;210)(1.15;225)
\psarc[liftpen=2](0,0){1.15}{315}{345}
\psline(1.15;315)(.8;330)(1.15;345)
\psbezier[liftpen=2](.4;90)(.3;130)(.3;170)(.4;210)
\psbezier(.4;210)(.3;250)(.3;290)(.4;330)
\psbezier(.4;330)(.3;10)(.3;50)(.4;90)}
\psline(1.15;75)(.8;90)(1.15;105)
\psline(1.15;195)(.8;210)(1.15;225)
\psline(1.15;315)(.8;330)(1.15;345)
\psline[doubleline=true](.4;90)(.8;90)
\psline[doubleline=true](.4;210)(.8;210)
\psline[doubleline=true](.4;330)(.8;330)
\psbezier(.4;90)(.3;130)(.3;170)(.4;210)
\psbezier(.4;210)(.3;250)(.3;290)(.4;330)
\psbezier(.4;330)(.3;10)(.3;50)(.4;90)
\endpspicture = 0
\label{etriang}
\end{equation}
To analyze it, we establish some conventions about graphs: In general, if $G$
is a graph, $V_G$ denotes the vertex set of $G$, and if $x,y \in V_G$, $xEy$ is
the relation that $x$ and $y$ are connected by an edge.  An injection $f:V_G
\to V_H$ is edge-respecting means that $f(x)Ef(y)$ if and only if $xEy$.  The
graph $G$ is $n$-point transitive means that if $H$ is a full subgraph of $G$
with at most $n$ vertices, every edge-respecting injection $f:V_H \to V_G$
extends to an automorphism of $G$.  More generally, $G$ is $n$-point regular
means the following:   For every $H$ with $k \le n$ vertices, for every $H'$
with $k+1$ vertices, for every edge-respecting map $f:V_H \to V_{H'}$, and for
every edge-respecting map $g:V_H \to V_G$, the number of ways to complete the
commutative diagram
$$
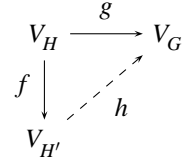
\begin{array}{cc}
\\
\Rnode{a}{V_H} & \hspace{1cm} \Rnode{b}{V_G} \\[1cm]
\Rnode{c}{V_{H'}}
\end{array}
\ncLine[nodesep=4pt,linecolor=black]{->}{a}{b}\Aput{g}
\ncLine[nodesep=4pt,linecolor=black]{->}{a}{c}\Bput{f}
\ncLine[nodesep=4pt,linecolor=black,linestyle=dashed,arrows=->]{c}{b}\Bput{h}
$$
with an edge-respecting map $h$ depends only on $H$, $H'$, and $f$ and not on
$g$.  Clearly, if $G$ is $n$-point transitive, then it is $n$-point regular.

Consider the numerical equalities implicit in equation~(\ref{etriang});
for any choice $(a,b,c)$ of three vertices of $J$, not necessarily
distinct, the left side becomes a state sum by labelling the
three outside regions by $a$, $b$, and $c$ and summing over the state
of the inside region:
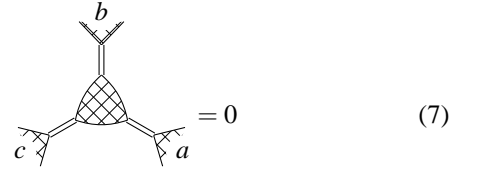
\begin{equation}
\pspicture[.4](-1.2,-1.2)(1.2,1.2)
\pscustom[fillstyle=crosshatch,linestyle=none]{
\psarc(0,0){1.15}{75}{105}
\psline(1.15;75)(.8;90)(1.15;105)
\psarc[liftpen=2](0,0){1.15}{195}{225}
\psline(1.15;195)(.8;210)(1.15;225)
\psarc[liftpen=2](0,0){1.15}{315}{345}
\psline(1.15;315)(.8;330)(1.15;345)
\psbezier[liftpen=2](.4;90)(.3;130)(.3;170)(.4;210)
\psbezier(.4;210)(.3;250)(.3;290)(.4;330)
\psbezier(.4;330)(.3;10)(.3;50)(.4;90)}
\psline(1.15;75)(.8;90)(1.15;105)
\psline(1.15;195)(.8;210)(1.15;225)
\psline(1.15;315)(.8;330)(1.15;345)
\psline[doubleline=true](.4;90)(.8;90)
\psline[doubleline=true](.4;210)(.8;210)
\psline[doubleline=true](.4;330)(.8;330)
\psbezier(.4;90)(.3;130)(.3;170)(.4;210)
\psbezier(.4;210)(.3;250)(.3;290)(.4;330)
\psbezier(.4;330)(.3;10)(.3;50)(.4;90)
\cput*(1.25;330){$a$}
\cput*(1.25;90){$b$}
\cput*(1.25;210){$c$}
\endpspicture = 0
\label{ethree}
\end{equation}
If $J$ is 3-point regular, this sum only depends on which of $a$, $b$ and $c$
are equal and which are connected by edges of $J$.

\begin{lemma}  If the graph $J$ of a checkerboard state model
is 3-point regular, equation~(\ref{etriang}) is a corollary of the
other checkerboard skein relations.
\end{lemma}

(Conversely, Jaeger~\cite{Jaeger:spin} proved that $J$ must be 3-point
regular if all of the checkerboard skein relations hold.)

\begin{proof} For convenience, we define a new type of strand,
denoted by dashes, as a linear combination of other webs:
$$
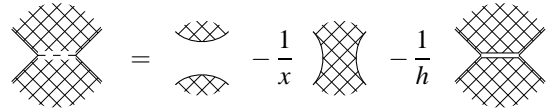
\pspicture[.4](-.9,-.7)(.9,.7)
\pscustom[fillstyle=crosshatch,linestyle=none]{
\psline(.6,-.35)(.25,0)(.6,.35)
\psarc(0,0){.695}{30.256}{149.744}
\psline(-.6,.35)(-.25,0)(-.6,-.35)
\psarc(0,0){.695}{-149.744}{-30.256}}
\psline[border=.05,linestyle=dashed](-.25,0)(.25,0)
\psline(-.6,-.35)(-.25,0)(-.6,.35)
\psline(.6,.35)(.25,0)(.6,-.35)
\endpspicture = \bwhh - \frac{1}{x}\bwvv - \frac{1}{h}\bwhhoriz
$$
This strand has its own weight matrix $W_D$ which is a linear combination of
$W_H$, $N$, and $I$; the weights are chosen so that $W_D(a,b) = 1$ if $a$ and
$b$ are distinct but not connected by an edge of $J$ and $W_D(a,b) = 0$ if
they are connected by an edge.  Then each case of equation~(\ref{ethree}) can
be converted to a statement about a state sum of a web on the sphere. For
example, if $a$, $b$, and $c$ form an anti-triangle and $J$ has $t$
anti-triangles, then the left side of equation~(\ref{ethree}) differs by a
factor of $t$ from the state sum of the graph
$$
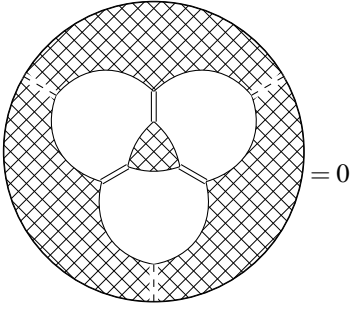
\pspicture[.4](-2,-2)(2,2)
\pscustom[fillstyle=crosshatch]{\evenodd
\psbezier(.4;90)(.3;130)(.3;170)(.4;210)
\psbezier(.4;210)(.3;250)(.3;290)(.4;330)
\psbezier(.4;330)(.3;10)(.3;50)(.4;90)
\psbezier[liftpen=2]( .8; 90)(1.3;110)(1.5;130)(1.5;150)
\psbezier(1.5;150)(1.5;170)(1.3;190)( .8;210)
\psbezier( .8;210)(1.3;230)(1.5;250)(1.5;270)
\psbezier(1.5;270)(1.5;290)(1.3;310)( .8;330)
\psbezier( .8;330)(1.3;350)(1.5; 10)(1.5; 30)
\psbezier(1.5; 30)(1.5; 50)(1.3; 70)( .8; 90)
\psarc[liftpen=2](0,0){2}{0}{360}}
\psline[doubleline=true](.4;90)(.8;90)
\psline[doubleline=true](.4;210)(.8;210)
\psline[doubleline=true](.4;330)(.8;330)
\psline[border=.05,linestyle=dashed](1.5;150)(2;150)
\psline[border=.05,linestyle=dashed](1.5;270)(2;270)
\psline[border=.05,linestyle=dashed](1.5; 30)(2; 30)
\pscircle(0,0){2}
\endpspicture = 0
$$
Let $w$ be either this web or its counterpart from one of the other cases of
equation~(\ref{ethree}).  Given that the $B_2$ skein relations are
consistent, and given that $w$ certainly does vanish modulo the skein
relations together, it suffices to show that $w$ is a multiple of the empty
web modulo all skein relations other than equation~(\ref{etriang}); the
coefficient is then necessarily zero.  Finally, since $w$ has at most four
black regions, it satisfies this condition by Lemma~\ref{lseven}.
\end{proof}

\begin{lemma} Any colored $B_2$ web $w$ on the sphere with at most seven black
regions is proportional to the empty web using only skein
relations other than equation~(\ref{etriang}). \label{lseven}
\end{lemma}
\begin{proof}
The proof is by induction on the number of vertices.  Following the usual
description of the $B_2$ spider \cite{Kuperberg:spiders}, we assign formal
angles of 45, 135 degrees, and 135 degrees to each vertex:
$$
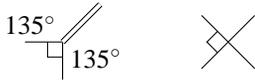
\pspicture(-.6,-.5)(.6,.5)
\psline(-.2,0)(-.2,-.2)(0,-.2)
\qline(-.5,0)(0,0)\qline(0,-.5)(0,0)
\psline[doubleline=true](0,0)(.5,.5)
\rput[rb](-.1,.1){$135^\circ$}
\rput[lt](.1,-.1){$135^\circ$}
\endpspicture
\hspace{1cm}
\pspicture(-.6,-.5)(.6,.5)
\psline(.2;135)(.2828;180)(.2;225)
\qline(.5;45)(.5;225)\qline(.5;135)(.5;315)
\endpspicture
$$
We use these formal angles to define a non-standard Euler characteristic of
each face of $w$ with the property that the total Euler characteristic is 2. 
(It should not be confused with the ordinary Euler characteristic used above
for normalization.)  It may be defined for a face as $1 - \frac{A}{360^\circ}$,
where $A$ is the sum of  all exterior angles, and an exterior angle at a given
vertex is defined as $180^\circ$ minus the usual interior angle. 
Alternatively, we may recall that the Euler characteristic of a face or a
vertex is 1 and that of an edge is $-1$, and then modify this definition by
dividing the Euler characteristic of an edge equally between the two faces that
contain it and, at each vertex, giving the right-angled faces $1/4$ of the
characteristic of the vertex and the other faces $3/8$ each.  This second
definition makes clear that the total is still 2.

Note also that the Euler characteristic of any face is a multiple of $1/4$. 
Therefore $w$ has either a white face with positive Euler characteristic or a
black face whose Euler characteristic is at least 1/2. Modulo webs with fewer
vertices, we can apply the exchange equation~(\ref{eexchange}) to all sides of
such a face which are type 2 strands, for example:
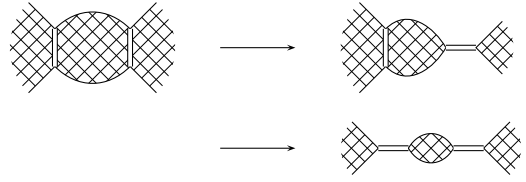
\begin{align*}
\pspicture[.5](-1.2,-.7)(1.2,.7)
\pscustom[fillstyle=crosshatch,linestyle=none]{
\psline(-.85,-.6)(-.5,-.25)
\psbezier(-.5,-.25)(-.2,-.55)(.2,-.55)(.5,-.25)
\psline(.5,-.25)(.85,-.6)
\psarc(.25,0){.849}{315}{45}
\psline(.85,.6)(.5,.25)
\psbezier(.5,.25)(.2,.55)(-.2,.55)(-.5,.25)
\psline(-.5,.25)(-.85,.6)
\psarc(-.25,0){.849}{135}{225}}
\psline(-.85,-.6)(-.5,-.25)
\psbezier(-.5,-.25)(-.2,-.55)(.2,-.55)(.5,-.25)
\psline(.5,-.25)(.85,-.6)
\psline(.85,.6)(.5,.25)
\psbezier(.5,.25)(.2,.55)(-.2,.55)(-.5,.25)
\psline(-.5,.25)(-.85,.6)
\psline[doubleline=true](-.5,-.25)(-.5,.25)
\psline[doubleline=true](.5,-.25)(.5,.25)
\endpspicture
&\psgoesto
\pspicture[.5](-1.2,-.7)(1.3,.7)
\pscustom[fillstyle=crosshatch,linestyle=none]{
\psarc(-.25,0){.849}{135}{225}
\psline(-.85,-.6)(-.5,-.25)
\psbezier(-.5,-.25)(-.2,-.55)(.15,-.26)(.3,0)
\psbezier(.3,0)(.15,.26)(-.2,.55)(-.5,.25)
\psline(-.5,.25)(-.85,.6)}
\psline(-.85,-.6)(-.5,-.25)
\psbezier(-.5,-.25)(-.2,-.55)(.15,-.26)(.3,0)
\psbezier(.3,0)(.15,.26)(-.2,.55)(-.5,.25)
\psline(-.5,.25)(-.85,.6)
\psline[doubleline=true](-.5,-.25)(-.5,.25)
\psline[doubleline=true](.3,0)(.7,0)
\pscustom[fillstyle=crosshatch,linestyle=none]{
\psarc(.7,0){.5}{315}{45}
\psline(1.05,.35)(.7,0)(1.05,-.35)}
\psline(1.05,.35)(.7,0)(1.05,-.35)
\endpspicture \\
&\psgoesto
\pspicture[.5](-1.3,-.5)(1.3,.5)
\pscustom[fillstyle=crosshatch,linestyle=none]{
\psarc(-.7,0){.5}{135}{225}
\psline(-1.05,-.35)(-.7,0)(-1.05,.35)}
\psline(-1.05,-.35)(-.7,0)(-1.05,.35)
\pscustom[fillstyle=crosshatch,linestyle=none]{
\psarc(.7,0){.5}{315}{45}
\psline(1.05,.35)(.7,0)(1.05,-.35)}
\psline(1.05,.35)(.7,0)(1.05,-.35)
\psline[doubleline=true](-.7,0)(-.3,0)
\pscustom[fillstyle=crosshatch]{
\psbezier(-.3,0)(-.15,.26)(.15,.26)(.3,0)
\psbezier(.3,0)(.15,-.26)(-.15,-.26)(-.3,0)}
\psline[doubleline=true](.3,0)(.7,0)
\endpspicture
\end{align*}
This operation does not change the Euler characteristic of the face. The face
can then be simplified by one of the other skein relations, for it has at
most three sides if it is white and at most two sides if it is black.
\end{proof}

There are only two known graphs $J$ that satisfy all of the above conditions,
namely the pentagon and the Higman-Sims graph.  The Higman-Sims graph is a
graph on $n=100$ vertices whose symmetry group contains the Higman-Sims
sporadic simple group as a subgroup of index two\cite{HS:group}. The
symmetry group acts 3-point transitively, so that the graph is 3-point
regular.  We claim that there is a corresponding state model that satisfies
the $B_2$ skein relations with $q^{1/2} = \tau$, which implies that $\ell =
-10$ and $h = -5$. Clearly $n = \ell^2$.  The graph has no triangles and each
vertex has valence $22 = \frac{n-x}{h}$, since $x = \ell$. The eigenvalues of
its adjacency matrix $A$ are 22, 2, and -8; the image of $N$ is the unique
line with eigenvalue 22. Therefore the eigenvalues of $W_I = hA + N - xI$ are
$-110+100+10 = 0$, $-10+10=0$, and $40+10 = 50$, which implies $W_I^2 =
xhW_I$, as desired.

\section{Discussion}

The properties of the Higman-Sims state model imply a number of mysterious
numerological connections between the Higman-Sims group ($HS$) and the (quantum)
representation theory of $\sp(4)$.  In this discussion, ``representation''
will in general mean a finite-dimensional linear representation.

The smallest representations of $\sp(4)$ are $V(\lambda_1)$, the defining
4-dimensional representation, $V(\lambda_2)$, the 5-dimensional representation
that identifies $\sp(4)$ with $\so(5)$, and $V(2\lambda_1)$, the 10-dimensional
representation which is the symmetric tensor square of $V(\lambda_1)$.  In
quantum representation theory, representations have a quantum dimension, which
is a natural generalization of the non-quantum or honest dimension. (It also
coincides with the character of a certain circle in the non-quantum
representation theory and appears in some proofs of the Weyl dimension formula
\cite{Bourbaki:lie}.) The quantum dimensions of these representations are
\begin{align*}
\dim_q V(\lambda_1)  &= q^2+q+q^{-1}+q^{-2} \\
\dim_q V(\lambda_2)  &= q^3 + q + 1 + q^{-1} + q^{-3} \\
\dim_q V(2\lambda_1) &= q^4+q^3+q^2+q+2+q^{-1}+q^{-2}+q^{-3}+q^{-4}
\end{align*}
At the same time, the quantum dimension is the value of a closed loop in the
$B_2$ spider; in the first two cases the loop is a type 1 or 2 strand, and in
the third case it is a dashed loop.  By the existence of the Higman-Sims state
model, these three numbers must also be 10, the square root of the number of
vertices of the Higman-Sims graph; 22, the degree of a vertex; and 77, the
anti-degree of a vertex.  Now the Higman-Sims graph has a special duality
realized by switching colors in the state model, and this duality tells us that
22 and 77 must also be dimensions of representations of $HS$; as it happens,
the two smallest irreducible representations.  Thus we learn that quantum
dimensions of representations of $U_{q=\tau^2}(\sp(4))$ coincide with honest
dimensions of representations of $HS$.  This pattern extends to all
representations of $U_{q=\tau^2}(\sp(4))$, except that the corresponding
representations of $HS$ are eventually not irreducible.

In fact, the $B_2$ spider can be understood as (the Hom spaces of) the
representation category of $U_q(\sp(4))$, and the Higman-Sims state model
establishes a functor from (the even half of) the representation category of
$U_{q=\tau^2}(\sp(4))$ to the representation category of $HS$.  Such a functor
would exist if there were an algebra homomorphism from the group algebra
$\C[HS]$ to the quantum group $U_{q=\tau^2}(\sp(4))$, but this possibility is
nonsense (as the referee mentioned), because it would relate honest dimensions
to honest dimensions and not quantum dimensions to honest dimensions.  But
perhaps one can construct $HS$ from $U_{q=\tau^2}(\sp(4))$ using this functor.

As a warm-up to this problem, one can try to construct a relationship between
the quantum group $U_q(\sl(2))$ and the symmetric group $S_n$ on $n = q + 2 +
q^{-1}$ letters, which is the symmetry group of the Potts model.  The Potts
model relates these objects in the same way that the Higman-Sims model relates
$U_{q=\tau^2}(\sp(4))$ and $HS$.

It would be especially interesting if one could construct not only the
Higman-Sims group but also other sporadic simple groups using quantum groups at
special values of $q$.

Besides its numerology, the Higman-Sims state model also has the following
interesting aspect. The most common axiomatic description of a spider, which is
a collection of web spaces for disks with different boundaries but with the
same skein relations, is that it is a kind of monoidal category, or a braided
category if crossings exist.  However, another interesting point of view is
that a spider is a certain kind of $2$-category with only one $0$-morphism (or
``object''), with a $1$-morphism for every choice of boundary, and such that
each web is a $2$-morphism \cite{BD:higher0}.  In this setting, a collection of
checkerboard skein modules has a suggestive definition also, namely as a
2-category with two 0-morphisms.

\acknowledgments

The author would like to thank Vaughan Jones, Pierre de la Harpe, and
Fran\c{c}ois Jaeger for fruitful discussions, as well as the reviewer and John
Baez for useful remarks and corrections.  The author used the \TeX\, macro
package PSTricks~\cite{pstricks} to typeset the equations and figures.

% \bibliography{qa,me,books,soft}

\begin{thebibliography}{10}

\bibitem{BD:higher0}
John~C. Baez and James Dolan, \emph{Higher-dimensional algebra and topological
  quantum field theory}, J. Math. Phys. \textbf{36} (1995), no.~11, 6073--6105,
  \mbox{arXiv:q-alg/9503002}.

\bibitem{Bourbaki:lie}
Nicolas Bourbaki, \emph{Groupes et alg\'ebres de {Lie}}, vol. 4--6, Hermann,
  Paris, 1968.

\bibitem{delaHarpe:jaeger}
Pierre de~la Harpe, \emph{Spin models for link polynomials, strongly regular
  graphs, and {Jaeger's} {Higman-Sims} model}, Pacific J. Math. \textbf{162}
  (1994), no.~1, 57--96.

\bibitem{dHJ:graph}
Pierre de~la Harpe and Vaughan F.~R. Jones, \emph{Graph invariants related to
  statistical mechanical models: examples and problems}, J. Combin. Theory Ser.
  B \textbf{57} (1993), no.~2, 207--227.

\bibitem{HS:group}
Donald~G. Higman and Charles~C. Sims, \emph{A simple group of order
  44,352,000}, Math. Zeitschr. \textbf{105} (1968), 110--113.

\bibitem{Jaeger:spin}
Fran\c{c}ois Jaeger, \emph{Strongly regular graphs and spin models for the
  {Kauffman} polynomial}, Geom. Dedicata \textbf{44} (1992), no.~1, 23--52.

\bibitem{Jones:pacific}
Vaughan F.~R. Jones, \emph{On knot invariants related to some statistical
  mechanical models}, Pacific J. Math. \textbf{137} (1989), no.~2, 311--334.

\bibitem{Kauffman:regular}
Louis~H. Kauffman, \emph{An invariant of regular isotopy}, Trans. Amer. Math.
  Soc. \textbf{318} (1990), no.~2, 417--471.

\bibitem{Kuperberg:g2}
Greg Kuperberg, \emph{The quantum {$G_2$} link invariant}, Internat. J. Math.
  \textbf{5} (1994), no.~1, 61--85.

\bibitem{Kuperberg:spiders}
\bysame, \emph{Spiders for rank 2 {Lie} algebras}, Comm. Math. Phys.
  \textbf{180} (1996), no.~1, 109--151, \mbox{arXiv:q-alg/9712003}.

\bibitem{pstricks}
\emph{{PSTricks}}, http://www.tug.org/applications/PSTricks/.

\end{thebibliography}

\providecommand{\bysame}{\leavevmode\hbox to3em{\hrulefill}\thinspace}

\end{document}